\documentclass[12pt,onecolumn]{IEEEtran}

\usepackage[colorlinks=false]{hyperref}
\usepackage[dvips]{graphicx}
\usepackage[usenames,dvipsnames]{color}
\usepackage{epsfig}
\usepackage{rotating}
\usepackage{amssymb}
\usepackage{amsmath,amsfonts}
\usepackage[boxed]{algorithm}
\usepackage[latin1]{inputenc}
\usepackage{verbatim}
\usepackage[tableposition=top,font=small,labelfont=bf,format=hang]{caption}
\usepackage{booktabs}
\newtheorem{lem}{Lemma}

\newtheorem{thm}{Theorem}
\newtheorem{prop}{Proposition}

\newtheorem{cor}{Corollary}

\title{ Robin inequality for $7-$free integers}
\author{   Patrick Sol\'{e}\thanks{Telecom ParisTech ,  46 rue Barrault, 75634
Paris Cedex 13, France.}\hspace{1cm} Michel Planat\thanks{ Institut FEMTO-ST, CNRS, 32 Avenue de
l'Observatoire, F-25044 Besan\c con, France.}}

\begin{document}
\maketitle

\begin{abstract} Recall that an integer is $t-$free iff it is not divisible by $p^t$ for some prime $p.$  We give a method to check Robin inequality $\sigma(n) < e^\gamma n\log\log n,$ for $t-$free integers $n$ and apply it for $t=6,7.$
We introduce $\Psi_t,$ a generalization of Dedekind $\Psi$ function defined for any integer $t\ge 2$ by
$$\Psi_t(n):=n\prod_{p \vert n}(1+1/p+\cdots+1/p^{t-1}).$$  
If $n$ is $t-$free then the sum of divisor function $\sigma(n)$ is $ \le \Psi_t(n).$ We characterize the champions for $x \mapsto \Psi_t(x)/x,$
as primorial numbers. Define the ratio $R_t(n):=\frac{\Psi_t(n)}{n\log\log n}.$ We prove that, for all $t$, there exists an integer $n_1(t),$ 
such that  we have $R_t(N_n)< e^\gamma$ for $n\ge n_1,$ where $N_n=\prod_{k=1}^np_k.$ Further, by combinatorial arguments, this can be extended to $R_t(N)\le e^\gamma$ for all
$N\ge N_n,$ such that  $n\ge n_1(t).$ This yields Robin inequality for $t=6,\,7.$ For $t$ varying slowly with $N$, we also derive $R_t(N)< e^\gamma.$

\end{abstract}

\section{Introduction}
The Riemann Hypothesis (RH), which describes the non trivial zeroes of Riemann $\zeta$ function has been deemed the Holy Grail of Mathematics by several authors \cite{Bo,L}.
There exist many equivalent formulations in the literature \cite{C}. The one of concern here is that of Robin \cite{R}, which is given in terms of $\sigma(n)$ the sum of divisor function

$$\sigma(n) < e^\gamma n\log\log n,$$
for $n\ge 5041.$ 
Recall that an integer is $t-$free iff it is not divisible by $p^t$ for some prime $p.$ 
The above inequality was checked for many infinite families of integers in \cite{CLMS}, for instance $5-$free integers.
In the present work we introduce  a method to check the inequality for $t-$free integers for larger values of $t$ and apply it to $t=6,7.$
The idea of our method is to introduce the generalized Dedekind $\Psi$ function defined for any integer $t\ge 2$ by
$$\Psi_t(n):=n\prod_{p \vert n}(1+1/p+\cdots+1/p^{t-1}).$$
If $t=2$ this is just the classical Dedekind function which occurs in the theory of modular forms \cite{modular}, in physics \cite{P}, and in analytic number theory \cite{SP}.
By construction, if $n$ is $t-$free then the sum of divisors function $\sigma(n)$
 is $ \le \Psi_t(n).$ To see this note that the multiplicative function $\sigma$ satisfies for any integer $a$ in the range $t>a\ge 2$
$$\sigma(p^{a})= 1+p+\cdots+p^a,$$
when the multiplicative function $\Psi_t$ satisfies
$$\Psi_t(p^{a})=p^a+\cdots+1+\cdots+1/p^{t-1-a}.$$
It turns out that the structure of champion numbers for the arithmetic function $x \mapsto \Psi_t(x)/x$ is much easier to understand than that
of $x \mapsto \sigma(x)/x,$ the super abundant numbers. They are exactly the so-called primorial numbers (product of first consecutive primes).
 We prove that, in order to maximize the ratio $R_t$ it is enough
to consider its value at primorial integers. Once this reduction is made, bounding above unconditionally $R_t$ is easy by using classical lemmas on partial eulerian products. 
We conclude the article by some results on $t-$free integers $N\ge N_n,$ valid for $t$ varying slowly with $N.$
\section{Reduction to primorial numbers}
Define the primorial number $N_n$ of index $n$ as the product of the first $n$ primes
$$N_n=\prod_{k=1}^np_k, $$
so that $N_0=1,N_1=2,\, N_2=6,\cdots$ and so on. 
The primorial numbers  (OEIS sequence $A002110$ \cite{OEIS}) play the role here of superabundant numbers in \cite{R} or primorials in \cite{N}.
They are champion numbers (ie left to right maxima) of the function $x \mapsto \Psi_t(x)/x:$
\begin{equation}
\frac{\Psi_t(m)}{m}<\frac{\Psi_t(n)}{n}~ \mbox{for}~\mbox{any}~  m<n. 
\label{newSuperab} 
\end{equation}
We give  a rigorous proof of this fact.
{\prop The primorial numbers and their multiples are exactly the champion numbers of the function $x \mapsto \Psi_t(x)/x.$ }
\begin{IEEEproof}
 The proof is by induction on $n$. The induction hypothesis $H_n$ is that the statement is true up to $N_n.$
Sloane sequence $A002110$ begins $1,2,4,6\dots$ so that $H_2$ is true. Assume $H_n$ true. Let $N_n\le m< N_{n+1}$ denote a generic integer. The prime divisors of $m$ are
$\le p_n.$ Therefore $\Psi_t(m)/m \le \Psi_t(N_n)/N_n$ with equality iff $m$ is a multiple of $N_n.$ Further $\Psi_t(N_n)/N_n< \Psi_t(N_{n+1})/N_{n+1}.$ The proof of 
$H_{n+1}$ follows.
\end{IEEEproof}

In this section we reduce the maximization of $R_t(n)$ over all integers $n$ to the maximization over primorials.

{\prop \label{reduc} Let $n$ be an integer $\ge 2.$ For any $m$ in the range  $N_n\le m< N_{n+1}$ one has $R_t(m)< R_t(N_n).$}
\begin{IEEEproof}
Like in the preceding proof we have 
$$\Psi_t(m)/m \le \Psi_t(N_n)/N_n.$$

Since $0<\log \log N_n \le \log \log m,$ the result follows.
\end{IEEEproof}

\section{$\Psi_t$ at primorial numbers}
We begin by an easy application of Mertens formula.

{\prop \label{mertens}  For $n$ going to $\infty$ we have $$\lim R_t(N_n)=\frac{e^\gamma}{\zeta(t)}.$$}

\begin{IEEEproof}
Writing $1+1/p=(1-1/p^2)/(1-1/p)$ in the definition of $\Psi(n)$
we can combine the Eulerian product for $\zeta(t)$ with Mertens formula
$$\prod_{p\le x}(1-1/p)^{-1}\sim e^\gamma \log(x)$$ to obtain

$$\Psi(N_n)\sim \frac{e^\gamma}{\zeta(t)} \log(p_n).$$
Now the Prime Number Theorem \cite[Th. 6, Th. 420]{HW} shows that $x \sim \theta(x)$ 
for $x$ large, where $\theta(x)$ stands for Chebyshev's first summatory function:
$$\theta(x)= \sum_{p\le x }\log p.$$ This shows
that, taking $x=p_n$ we have
$$p_n\sim \theta(p_n)=\log(N_n).$$
The result follows.
\end{IEEEproof}

This motivates the search for explicit upper bounds on $R_t(N_n)$ of the form $\frac{e^\gamma}{\zeta(t)}(1+o(1)).$
In that direction we have the following bound.

{\prop \label{fonda} For $n$ large enough to have $p_n\ge 20000,$ we have $$ \frac{\Psi_t(N_n)}{N_n} \le \frac{\exp(\gamma+2/p_n)}{\zeta(t)}(\log\log N_n+\frac{1.1253}{\log p_n}) .$$}
We prepare for the proof of the preceding Proposition by some Lemmas.
First an upper bound on a partial Eulerian product from \cite[(3.30) p.70]{RS}.
{\lem \label{rs} For $x\ge 2,$ we have $$\prod_{p \le x } (1-1/p)^{-1}\le e^{\gamma} (\log x+\frac{1}{\log x}).$$}
Next an upper bound on the tail of the Eulerian product for $\zeta(t).$
{\lem \label{us} For $n\ge 2$ we have  $$\prod_{p >p_n } (1-1/p^t)^{-1}\le \exp(2/p_n).$$}

\begin{IEEEproof}
Use Lemma 6.4 in \cite{CLMS} with $x=p_n.$ Bound $\frac{t}{t-1}x ^{1-t}$ above by $2/x.$
\end{IEEEproof}

{\lem \label{robmod} For $n\ge 2263,$ we have

$$\log p_n<\log \log N_n+ \frac{0.1253}{\log p_n}. $$}

\begin{IEEEproof}
If $n\ge 2263,$  then $p_n\ge 20000.$ By \cite{RS}, we know then that

$$ \log N_n >p_n(1-\frac{1}{8p_n}).$$
 On taking log's we obtain 
$$ \log\log N_n >\log p_n-\frac{0.1253}{p_n},$$
upon using
$$\log(1-\frac{x}{8}) >  -0.1253 x $$
for $x$ small enough. In particular $x<1/20000$ is enough.
\end{IEEEproof}

We are now ready for the proof of Proposition \ref{fonda}.
\begin{IEEEproof}

Write $$\frac{\Psi_t(N_n)}{N_n} = \prod_{k=1}^n \frac{1-1/{p_k}^t}{1-1/p_k}=\frac{\prod_{p >p_n } (1-1/p^t)^{-1}}{\zeta(t)}\prod_{p \le p_n} (1-1/p)^{-1}$$ and use both Lemmas to derive

$$ \frac{\Psi_t(N_n)}{N_n} \le \frac{\exp(\gamma+2/p_n)}{\zeta(t)}(\log p_n+\frac{1}{\log p_n}). $$

Now we get rid of the first $\log$ in the RHS by Lemma \ref{robmod}.

The result follows.
 
\end{IEEEproof}

So, armed with this powerful tool, we derive the following significant Corollaries.\\
For convenience let 
$$ f(n)=(1+\frac{1.1253}{\log p_n \log\log N_n}).$$
{\cor \label{fonda+} Let $n_0=2263.$  Let $n_1(t)$ denote the least $n\ge n_0$ such that $e^{2/p_n}f(n)<\zeta(t).$ 
For $n\ge n_1(t)$ we have $R_t(N_n)< e^\gamma.$ }

\begin{IEEEproof}

Let $n\ge n_0.$ We need to check that 
$$ \exp(2/p_n) (1+\frac{1.1253}{\log p_n  \log\log N_n})\le \zeta(t). $$
which, for fixed $t$  holds for $n$ large enough. Indeed $\zeta(t)>1$ and the LHS goes monotonically to $1^+$ for $n$ large.

\end{IEEEproof}

We give a numerical illustration of Corollary \ref{fonda+} in Table 1.
\begin{table}[ht]
\begin{center}
\small
\begin{tabular}{|r|r|r|}
\hline
$t$ & $n_1(t)$ & $N_{n_1(t)}$ \\
\hline
3 & 10& $6.5 \times 10^{9}$\\
4 & 24& $2.4\times10^{34}$\\
5& 79&$4.1 \times 10^{163}$ \\
6& 509& $ 5.8 \times 10^{1551}$\\
7&10 596& $2.5\times 10^{48337}$\\

\hline
\end{tabular}
\label{table1}
\caption{The numbers in Corollary \ref{fonda+}.}
\end{center}
\end{table}

We can extend this Corollary to all integers $\ge n_0$ by using the reduction of preceding section.

{\cor \label{upper} For all $N\ge N_n$ such that $n\ge n_1(t)$ we have $R_t(N)< e^\gamma.$ }\\

\begin{IEEEproof}
 Combine Corollary \ref{fonda+} with Proposition \ref{reduc}.
\end{IEEEproof}

We are now in a position to derive the main result of this note.

{\thm If $N$ is a $7-$free integer, then $\sigma(N)<N e^\gamma \log \log N.$}

\begin{IEEEproof}
 If $N$ is $\ge N_n$ with $n\ge n_1(7),$ then the above upper bound holds for $\Psi_7(N)$ by Corollary \ref{upper}, hence for $\sigma(N)$ by the remark in the Introduction.
 If not, we invoke the results of \cite{B}, who checked Robin inequality for $5040<N\le 10^{10^{10}},$ and observe that all $7-$free integers are $>5040.$
\end{IEEEproof}

\section{Varying $t$}
We begin with an easy Lemma.
{\lem \label{geo} Let $t$ be a real variable. For  $t$ large, we have $\zeta(t)=1+\frac{1}{2^t}+o(\frac{1}{2^t}).$}
\begin{IEEEproof}
By definition, for $t>1$ we may write $$\zeta(t)=\sum_{n=1}^\infty\frac{1}{n^t}$$ so that
$$\zeta(t)\ge 1+\frac{1}{2^t}.$$
In the other direction, we write
$$\zeta(t)=1+\frac{1}{2^t}+\frac{1}{3^t}+\sum_{n=4}^\infty\frac{1}{n^t},$$
and compare the remainder of the series expansion of the $\zeta$ function with an integral:
$$ \sum_{n=4}^\infty\frac{1}{n^t}<\int_3^{\infty}\frac{du}{u^t}=\frac{3}{(t-1)3^t}=O( \frac{1}{3^t}).$$
The result follows.
\end{IEEEproof}

We can derive a result when $t$ grows slowly with $n.$
{\thm \label{S} Let $S_n$ be a sequence of integers such that $S_n\ge N_n$ for $n$ large, and such that $S_n$ is $t-$free with $t=o(\log \log n).$ For $n$ large enough, 
Robin inequality holds for $S_n.$}
\begin{IEEEproof}
 For Corollary \ref{upper} to hold we need 
 $$ e^{2/p_n}f(n)< \zeta(t)$$ to hold, or , taking logs, the exact bound
$$ 2/p_n+\log f(n) <\log \zeta(t), $$ 
or up to $o(1)$ terms

$$2/p_n+\frac{1.1253}{\log p_n \log\log N_n}\le \log \zeta(t).$$

In the LHS, the dominant term is of order $1/(\log p_n)^2,$ since, like in the proof of Proposition \ref{mertens}, we may write  $p_n\sim \log N_n$ . Now $p_n\sim n\log n$ by \cite[Th. 8]{HW}, entailing
$\log p_n\sim \log n$ and $(\log p_n)^2\sim (\log n)^2$. In the RHS, with the hypothesis made on $t$ we have, by Lemma \ref{geo}, the estimate $\log \zeta(t)\sim \frac{1}{2^t}.$ The result follows after comparing logarithms
of both sides.
\end{IEEEproof}
\section{Conclusion}

In this article we have proposed a technique to check Robin inequality for $t-$free integers for some values of $t.$ The main idea has been to investigate the complex structure of the divisor function $\sigma$ though the sequence of Dedekind psi functions $\psi_t$. The latter are simpler for the following reasons

\begin{itemize}
 \item $\Psi_t(n)$ solely depends on the prime divisors of $n$ and not on their multiplicity
\item the champions of $\Psi_t$ are the primorials instead of the colossally abundant numbers
\item $\Psi_t$ is easier to bound for $n$ large because of connections with Eulerian products
\end{itemize}

Further, $\sigma(n)\le \Psi_t(n)$ for $t-$free integers $n.$ We checked Robin inequality for $t-$free integers  for $t=6,7$ and $t=o(\log \log n).$ It is an interesting and difficult open problem
to apply Theorem \ref{S} to superabundant numbers or colossally abundant numbers for instance. We do not believe it is possible. New ideas are required to prove Robin inequality in full generality.


\end{document}